\documentclass[a4paper, 12pt,oneside,reqno]{amsart}
\usepackage[a4paper]{geometry}
\usepackage[matrix,arrow,curve,cmtip]{xy}

\usepackage{txfonts}
\DeclareMathAlphabet{\mathcal}{OMS}{cmsy}{m}{n} 

\usepackage{
  amsmath,  
  amssymb,  
  amsthm,   
  tikz,     
  fullpage, 
  youngtab, 
  ytableau, 
  thmtools, 
  hyperref, 
  cite,     
  verbatim,
  url       
}

\usepackage{amscd}
\usepackage{euscript}
\usepackage{latexsym}
\usepackage[cp1251]{inputenc}
\usepackage[english]{babel}
\usepackage{mathrsfs}
\usepackage{graphicx}
\usepackage{keyval}
\usepackage{mathtools}
\usepackage{tikz}
\usetikzlibrary{arrows,shapes,snakes,automata,backgrounds,petri,through,positioning}
\usetikzlibrary{intersections}
\usepackage[symbol*]{footmisc}
\usepackage{verbatim, xcolor}
\usepackage{tikz-cd}

\DeclareMathOperator{\A}{A}
\DeclareMathOperator{\B}{B}
\DeclareMathOperator{\C}{C}

\DeclareMathOperator{\dcobound}{d}

\DeclareMathOperator{\Hom}{Hom}
\DeclareMathOperator{\Homol}{H}

\def\GSB{\mathrm{GSB}}

\usepackage{etoolbox}
\patchcmd{\thebibliography}{\section*}{\paragraph}{}{}

\numberwithin{equation}{section}

\def\Vir{\mathrm{Vir}}
\def\Der{\mathop {\fam 0 Der} \nolimits}

\def\oo#1{\mathbin {{}_{(#1)}}}

\newtheorem{theorem}{Theorem}[section]
\newtheorem{lemma}[theorem]{Lemma}
\newtheorem{proposition}[theorem]{Proposition}
\newtheorem{corollary}[theorem]{Corollary}

\theoremstyle{definition}
\newtheorem{definition}[theorem]{Definition}
\newtheorem{example}[theorem]{Example}

\title{Hochschild cohomology of the Weyl conformal algebra with coefficients in finite modules}
\author{H. Alhussein$^{1),2)}$}
\author{P. Kolesnikov$^{3)}$}

\address{$^{1)}$ Siberian State University of Telecommunication and Informatics, Novosibirsk, Russia.}
\address{$^{2)}$Novosibirsk State University of Economics and Management, Novosibirsk, Russia.}
\address{$^{3)}$Sobolev Institute of Mathematics, Novosibirsk, Russia.}

\def\Vir{\mathop {\fam 0 Vir} \nolimits}
\def\Cend{\mathop {\fam 0 Cend} \nolimits}
\def\oo#1{\mathbin {{}_{(#1)}} }

\begin{document}
\maketitle
\begin{abstract}
In this work we find Hochschild cohomology groups of the 
Weyl associative conformal algebra with coefficients in all finite modules. 
The Weyl conformal algebra is the universal 
associative conformal envelope of the Virasoro Lie conformal algebra 
relative to the locality $N=2$. 
In order to obtain this result we adjust the  
algebraic discrete Morse theory 
to the case of differential algebras.
\end{abstract}

\section{Introduction}

Conformal algebras were introduced by Kac \cite{KacValgBeginners}
to formalize the properties of the coefficients of the 
singular part of the operator product expansion  (OPE) 
 in conformal field theory. 
In particular, every vertex algebra \cite{bor} 
is actually a Lie conformal algebra.

Assume $V$ is a vertex algebra with a translation operator 
$\partial$ and a state-field correspondence $Y$. 
Then, as a result of the locality axiom, the operator product expansion of two fields $Y (a, z)$ and $Y(b, z)$, $a, b \in V$, has a finite singular part:
\[
Y (a, w)Y (b, z) =\sum^{N(a,b)-1}_{n=0} Y (c_n, z)\frac{1}{(w-z)^{n+1}}+ (\text{regular part}).
\]
The coefficients of the singular part are determined by the commutator of the fields:
\[
[Y (a, w),Y (b, z)] =\sum^{N(a,b)-1}_{n=0} Y (c_n, z)\frac{1}{n!}\frac{\partial^n\delta(w-z)}{\partial z^n}.
\]
where $\delta(w-z)=\sum_{s\in Z} w^sz^{-s-1} $is the formal delta-function. The correspondence
\[
(a,b) \mapsto c_n \ , \ n\geq 0.
\]
defines an infinite series of bilinear operations ($n$-products) on $V$. Together with the translation operator $\partial$, these operations turn $V$ into a system called {\em conformal Lie algebra}. The structure theory of finite Lie conformal algebras and superalgebras was developed in \cite{DK1998, Fattori}. 
Cohomology theory of conformal algebras was introduced in \cite{BKV}, then developed in a more general context of pseudo-algebras in \cite{BDK}. 

The study of universal structures for conformal algebras was initiated in \cite{Roitman99}. 
The classical theory of Lie and associative algebras often needs universal constructions like free algebras and universal enveloping algebras. This was a motivation 
for the development of combinatorial issues in the theory of conformal algebras, in particular, Gr\"obner--Shirshov bases 
theory \cite{BFK}.

For every conformal Lie algebra $L$ 
one can construct a series of universal enveloping 
associative conformal algebras corresponding to 
different associative locality functions on the 
generators \cite{Roit2000}. 
For example, consider the Virasoro conformal algebra $\Vir $
which is generated by a single element~$v$.
One may fix a natural number $N$ 
and construct 
the associative conformal algebra $U(N)$ 
generated by the element $v$ such that 
$(v\oo{n} v) = 0$ for $n\ge N$, and 
the commutation relations of $\Vir $ hold.
Obviously, $U(1)=0$; the algebra $U(2)$ 
is known as the Weyl conformal algebra 
(also denoted $\Cend_{1,x}$ in \cite{BKL}). 

For every conformal algebra $C$ one may construct 
an ``ordinary'' algebra $\mathcal A(C)$ called a 
coefficient algebra of~$C$, which inherits many properties of $C$ (associative, commutative, Lie, Jordan, etc.).
Every conformal module over $C$ is also a module over 
$\mathcal A(C)$. Moreover, it was proved in \cite{BKV} that the (reduced) cohomology of a conformal algebra $C$ may be calculated via the corresponding cochain complex of its coefficient $\mathcal A(C)$.

It was shown in \cite{Kozlov2017} that the second Hochschild 
cohomology groups $\Homol ^2(U(2), M)$ are trivial for every conformal 
(bi-)module~$M$, but for 
higher Hochschild cohomologies the direct computation becomes too complicated. 
In contrast to the ``ordinary'' Hochschild cohomology, 
if $C$ is an infinite associative conformal algebra 
then 
one cannot reduce the computation of 
$\Homol^n(C,M)$ to $\Homol^{n-1}(C, \mathrm{Chom}(C,M) )$
since the space of conformal homomorphisms 
$\mathrm {Chom}(C,M)$ may not be a conformal module 
over~$C$. Even if $M$ is finite (in this case, 
 the conformal $C$-module 
 $\mathrm {Chom}(C,M)$ is infinite and thus we may just
derive from the result of \cite{Kozlov2017}
that $\Homol^3(U(2),M)=0$ for every finite 
module~$M$.

In this paper we find all higher Hochschild 
cohomology groups $\Homol^n(U(2),M)$, $n\ge 2$,
of the Weyl conformal 
algebra $U(2)$ with 
coefficients in all finite modules~$M$. In order to obtain this result we construct the Anick resolution for its coefficient algebra via the algebraic discrete Morse theory as presented, for example, 
in \cite{JW}.
It is discussed in \cite{Akl} how to adjust this technique for differential algebras to calculate Hochschild cohomologies with coefficients in a trivial module. 
The purpose of our work is to apply the Morse matching method
for calculation of Hochschild cohomologies
of associative conformal algebras with coefficients 
in an arbitrary module.
As a result, we find that all Hochschild cohomology 
groups 
of $U(2)$ with coefficients in a finite module 
are trivial except the first one.
This result is close to the observation of \cite{SuY} 
about cohomologies of the Lie conformal algebra 
$\mathrm{gc}_n$ with coefficients in the natural 
module. 

\section{Preliminaries in conformal algebras}

Throughout the paper, 
$\Bbbk$ is a field of characteristic zero,
$\mathbb Z_+$ is the set of nonnegative integers.

\begin{definition} A conformal algebra \cite{KacValgBeginners} is a linear space $C$ equipped with a linear map 
$\partial:C \rightarrow C$ 
and a family of bilinear operations 
$(\cdot \oo n \cdot) :C\otimes C \rightarrow C$, 
$n \in \mathbb Z_+ $,
satisfying the following  properties:

(C1) for every $a,b \in C$ there exists $N = N(a,b)\in \mathbb Z_+$ such that $(a \oo n b) = 0$ for all $n\geq N$ ;\\
(C2) $(\partial a\oo n b) = -n(a \oo {n-1} b)$;\\
(C3)  $(a \oo n \partial b) = \partial(a \oo n b) + n(a \oo {n-1} b)$.

 Every conformal algebra $C$ is a left module over the polynomial algebra $H=\Bbbk [\partial ]$.
The structure of a conformal algebra on an $H$-module $C$
may be expressed by means of
a single polynomial-valued map called $\lambda$-product: 
\[
\begin{gathered}
(\cdot \oo{\lambda} \cdot) : C \otimes C \longrightarrow C[\lambda],
\\
(a \oo{\lambda} b) =\sum_{n=0}^{N(a,b)-1} \frac{\lambda^{n}}{n!}(a \oo{n} b),
\end{gathered}
\] 
where $\lambda $ is a formal variable, 
satisfying the following axioms:
\begin{gather}
    (\partial a\oo\lambda b) = -\lambda (a\oo\lambda b), 
      \label{eq:3/2-lin(1)}\\
    (a\oo\lambda \partial b) = (\partial+\lambda) (a\oo\lambda b).
      \label{eq:3/2-lin(2)}
\end{gather}
The number $N=N(a,b)$ is said to be a {\em locality level} of $a,b\in C$.
\end{definition}

\begin{definition}
 For every conformal algebra $C$ one may construct 
an algebra $A=\mathcal A(C)$ in the following way. 
As a linear space,
\[
A = \Bbbk [t,t^{-1}] \otimes _{\Bbbk [\partial ]} C,
\]
where $\partial $ acts on $\Bbbk[t,t^{-1}]$ as $- d/dt$.
Denote $a(n)=t^n\otimes_{\Bbbk [\partial ]} a$ for $a\in C$, $n\in \mathbb Z$.
The operation on $A$ is given by a well-defined expression
\[
a(n)\cdot b(m) = \sum\limits_{s\ge 0} \binom{n}{s} (a\oo{s} b) (n+m-s).
\]
The algebra $\mathcal A(C)$ is called the {\em coefficient algebra} of~$C$.
\end{definition}

The coefficient algebra $\mathcal A(C)$ has a derivation 
also denoted $\partial $ such that 
$\partial (a(n)) = (\partial a)(n) = -na(n-1)$.
The space $\mathcal A_+(C)$ spanned by all $a(n)$, $n\in \mathbb Z_+$, $a\in C$,
is a subalgebra of $\mathcal A(C)$ which is closed under the derivation~$\partial $.

Conformal algebra $C$ is called associative 
(commutative, Lie, Jordan, etc.) if so is $A(C)$ \cite{Roitman99}. 
For example, $A(C)$ is associative if and only if 
\begin{equation}\label{eq:ConfAss}
a \oo n (b \oo m c) = \sum _{s \geq 0} \binom{n}{s}
(a \oo{n-s} b) \oo {m+s} c
\end{equation}
for all $a, b, c \in C, n,m \in \mathbb Z_+$. 
In terms of the $\lambda$-product, the last relation may be expressed by a single formula
\[
a \oo \lambda (b \oo \mu c) = (a \oo \lambda b) \oo {\lambda+\mu} c, \quad a, b, c \in C,
\]
where $\lambda$ and $\mu$ are independent commuting variables \cite{KacForDistrib}.

In general, an associative conformal algebra $C$ turns into a Lie conformal algebra $C^{(-)}$
when equipped with a new $\lambda $-product $[\cdot\oo\lambda \cdot ]$ defined as follows:
\[
 [a\oo\lambda b] = (a\oo\lambda b) - (b\oo{-\partial - \lambda } a),\quad a,b\in C.
\]

\begin{example}
(1) The 1-generated free $H$-module $C = \Bbbk [\partial ] v$
equipped with 
\[
 (v\oo\lambda v ) = (\partial +2\lambda )v
\]
is a Lie conformal algebra known as the Virasoro conformal algebra.
The coefficient algebra of $\Vir $ is the Witt algebra
 $\mathcal A(\Vir) = \Der \Bbbk [t,t^{-1}]$, 
 and $\mathcal A_+(\Vir )= \Der \Bbbk [t]$.

(2) The infinitely generated free $H$-module 
$C = 
H\otimes \Bbbk\{v,v^2,v^3,\dots  \} \simeq \Bbbk [\partial, v]v$
equipped with 
\[
 (v^n\oo\lambda v^m) = v^n(v+\lambda )^m, \quad n,m\ge 1
\]
is an associative conformal algebra denoted ${U(2)}$.
The coefficient algebra $\mathcal A (U(2))$ 
is known to be the right ideal in the 
associative algebra $\Bbbk \langle p,q,q^{-1} \mid qp-pq=1\rangle $
generated by $p$. The positive part $\mathcal A_+(U(2))$
is isomorphic to $pW_1$, where $W_1$ is the first Weyl algebra 
generated by $p$, $q$ with $[q,p]=1$.
This is why $U(2)$ is called the {\em conformal Weyl algebra} \cite{Roitman99}.
\end{example}

Note that $\Vir \subset U(2)^{(-)}$, and $U(2)$ is generated 
(as an associative conformal algebra) by the elements of $\Vir $.
Hence, $U(2)$ is an associative envelope of $\Vir $. 
Moreover, this envelope is universal in the class 
of all associative envelopes $C$ of $\Vir $ such that 
$N(v,v)\le 2$ in $C$ \cite{BFK}. 
In fact, for the Virasoro conformal algebra $\Vir $ one may construct a series of associative conformal algebras $U(N)$, $N\in \mathbb Z_+$. 
Each $U(N)$ is the universal enveloping associative conformal algebra of $\Vir $ in the class of all associative conformal envelopes $C$
of $\Vir $ with $N(v,v)\le N$.

\begin{definition}[\cite{ChengKac}] A module $M$ over an associative conformal algebra $C$ is a $\Bbbk[\partial]$-module endowed with the $\lambda$-action $a_{(\lambda)}m$ which defines a map 
$C \otimes M\rightarrow M [\lambda] $ such that:
\begin{gather}
    (\partial a\oo\lambda m) = -\lambda (a\oo\lambda m), 
      \quad 
    (a\oo\lambda \partial m) = (\partial+\lambda) (a\oo\lambda m),
      \label{eq:3/2-lin(mod1-2)}\\
    (a\oo{\lambda } (b\oo\mu m)) = ((a\oo{\lambda } b)\oo{\mu+\lambda } m).
\end{gather}

Similarly, a conformal action of a Lie conformal algebra $L$
on a module $M$ meets \eqref{eq:3/2-lin(mod1-2)} 
and the conformal analogue of the Jacobi identity:
\[
(a\oo\lambda (b\oo\mu m )) - (b\oo\mu (a\oo\lambda m)) =
((a\oo\lambda b)\oo{\lambda+\mu} m),
\]
for $a,b\in L$, $m\in M$.
\end{definition}

\begin{example}
Given a 1-generated free $H$-module  $M = \Bbbk [\partial ]u$
and two scalars $\Delta, \alpha \in \Bbbk$, 
one may define conformal action of $\Vir $  on $M$ as
\begin{equation}\label{eq:Vir}
(v\oo\lambda u) = (\alpha +\partial +\Delta\lambda )u.
\end{equation}
Denote the conformal $\Vir$-module obtained by $M_{(\alpha, \Delta)}$.
For $\Delta \ne 0$, this is an irreducible $\Vir $-module, 
and every finite irreducible $\Vir $-module is isomorphic 
to an appropriate $M_{(\alpha, \Delta)}$ \cite{ChengKac}.
\end{example}

If $M$ is a module over an (associative or Lie) conformal algebra 
$C$ then $M$ is also a module over the ordinary 
(associative or Lie, respectively)
algebra $\mathcal A_+(C)$. 
Namely, for $a\in C$, $n\in \mathbb Z_+$, $u\in M$ the element 
$a(n)u$ is the coefficient at $\lambda ^{n}/n!$ of $(a\oo\lambda u)$:
\[
a\oo\lambda u  = \sum\limits_{n\ge 0} \dfrac{\lambda^n}{n!} a(n)u.
\]

For every conformal $C$-module $M$ over an associative conformal algebra $C$, the space $M$ is also a $C^{(-)}$-module relative to the same conformal action. 
The converse construction has a restriction due to locality. 
For example, the module $M_{(\alpha,\Delta)}$
over the Virasoro (Lie) conformal algebra is also a module 
over its universal enveloping associative conformal algebra 
$U(2)$
if and only if $\Delta = 0$ or $\Delta =1$.

Indeed, $(v\oo\lambda v) = v^2+\lambda v$ in $U(2)$, so 
\[
(v\oo\lambda v)\oo\mu u = v\oo\lambda (v\oo{\mu-\lambda } u )
=
v\oo\lambda (\alpha +\partial + \Delta(\mu-\lambda ))u
=(\alpha +\partial +\lambda + \Delta(\mu-\lambda ))(\alpha+\partial +\Delta\lambda )u.
\]
The polynomial in the right-hand side is of degree $<2$ 
in $\lambda $ if and only 
if $\Delta=0$ or $\Delta =1$. 

In \cite{BKL}, the algebra $U(2)$ (up to an isomorphism) is denoted 
by $\Cend_{1,x}$. The classification of finite irreducible modules 
over $\Cend_{1,x}$ indeed consists of the modules $M_{(\alpha, 1)}$.
The modules $M_{(\alpha,0)}$ are not irreducible, they contain 
submodules $(\partial+\alpha)M_{(\alpha, 0)}$.

The explicit formula for the action of $vh(\partial, v)\in U(2)$ on 
$M_{(a,1)}=Hu$, $H=\Bbbk [\partial ]$,
is a consequence of \eqref{eq:Vir} and associativity \eqref{eq:ConfAss}:
\[
vh(\partial, v)\oo\lambda f(\partial )u = 
h(-\lambda , \partial+\alpha )(\lambda +\partial +\alpha )
f(\partial +\lambda )u, 
\quad f(\partial )\in H.
\]
For example, the element $v(n)f(\partial)u$, $n\in \mathbb Z_+$, 
$f(\partial )\in H$,  is given by
\begin{equation}\label{eq:v(n)-Action}
v(n) f(\partial ) u  = (\partial+\alpha ) f^{(n)}(\partial )u+n f^{(n-1)}(\partial )u,
\end{equation}
where $f^{(n)}(\partial )$ 
stands for the $n$th derivative of $f(\partial )$.

Even for finite conformal modules over a simple associative conformal algebra,
there is no complete reducibility. 
For example, if $M_{(\alpha, 1)}$ and $M_{(\beta, 1)}$ are two irreducible 
modules over the Weyl conformal algebra $U(2)$ then 
there exists a non-split extension 
\[
0\to M_{(\alpha, 1)}=Hu \to E \to M_{(\beta , 1)}=Hw \to 0,
\]
e.g., 
\[
(v\oo\lambda w) = (\lambda +\partial +\beta)w + \gamma u,
\]
for $\gamma \in \Bbbk $.
Indeed, all extensions like that for the Virasoro (Lie) conformal algebra 
were described in \cite{ChengKac}, it is enough to choose those 
of limited locality.

\subsection{Hochschild cohomology for associative conformal algebras}\label{2.2}

Let $C$ be an associative conformal algebra, 
and let $M$ be a conformal module over~$C$.
The {\em basic Hochschild complex} \cite{BKV} 
$\tilde \C^\bullet (C,M)$
consists of the cochain spaces $\tilde \C^n (C,M)$, $n=1,2,\ldots$,
each of them is the space of all maps 
\[
\varphi_{\bar \lambda }: 
C^{\otimes n}\to M[\bar \lambda ],
\]
where $\bar\lambda = (\lambda_1,\ldots, \lambda_{n})$,
satisfying the conformal anti-linearity condition: 
\[
 \varphi_{\bar \lambda }(a_1,\ldots,\partial a_i, \ldots , a_n) =
 -\lambda_i \varphi _{\bar \lambda }(a_1,\ldots, a_n),\quad i=1,\ldots, n.
\]
The Hochschild differential 
$\dcobound_n : \tilde \C^n(C,M) \to \tilde \C^{n+1}(C,M)$
on the basic complex is given by
\begin{multline}\nonumber
(\dcobound_n\varphi)_{\bar\lambda }
(a_1,\ldots, a_{n+1}) = 
a_1 \oo{\lambda_1} 
\varphi_{\bar\lambda_0} (a_2, \ldots , a_{n+1})  
 + \sum\limits_{i = 1}^{n} (-1)^i\varphi_{\bar\lambda_i}
  (a_1,\ldots , a_i \oo{\lambda_i} a_{i+1},\ldots , a_{n+1}),
\end{multline}
for 
$\bar \lambda = (\lambda_1, \ldots, \lambda_{n+1})$, 
$\bar\lambda_0 = (\lambda_2, \ldots, \lambda_{n+1})$, 
$\bar\lambda_i = (\lambda_1, \ldots, \lambda_i+\lambda_{i+1},
\ldots , \lambda_{n+1} )$, $i=1,\dots, n$. The cohomology of the  basic Hochschild complex is called {\em the basic Hochschild cohomology} $\tilde H^\bullet (C, M)$. 

For every $n\ge 1$, the cochain space $\tilde \C^n(C,M)$ 
is a left  $\Bbbk[\partial]$-module:
\[
(\partial\varphi)_{\bar \lambda}(a_1,\ldots,a_n)=(\partial+\sum_{i=1}^{n}\lambda_i)\varphi_{\bar \lambda}(a_1,\ldots,a_n).
\]
For every $n\ge 1$, the map $\dcobound_n$ commutes with $\partial$. 
The quotient complex
\[
\C^\bullet (C,M)= 
\tilde \C^\bullet (C, M)/\partial \tilde \C^\bullet (C,M)\]
is called the {\em reduced Hochschild complex} and its cohomology is called the {\em reduced  Hochschild cohomology} $H^\bullet (C, M)$.

Another approach to the definition of Hochschild cohomologies for 
associative conformal algebras was considered in \cite{BDK}, 
see also \cite{Dolg2007}. The Hochschild cohomology groups 
defined there coincide with the reduced Hochschild cohomologies 
when written in terms of $\lambda $-operations.

Consider the ``ordinary'' Hochschild complex 
$\C^\bullet (\mathcal A_+(C),M)$
The maps
\[
\partial_n^* : \C^n( \mathcal A_+(C),M) \to \C^n( \mathcal A_+(C) ,M)
\]
given by
\[
(\partial_n^* f)(\alpha_1,\ldots, \alpha_n)
= \partial f(\alpha_1,\ldots, \alpha_n)
- \sum\limits_{i=1}^n 
f(\alpha_1,\ldots, \partial\alpha_i, \ldots, \alpha_n),
\]
and $\partial a(n) = -n a(n-1)$ for $a\in C$, $n\ge 0$, 
commute with the Hochschild differentials.

The following statement describes the relations between 
Hochschild cohomologies of $C$ and $\mathcal A_+(C)$ 
with coefficients on the same module~$M$.

\begin{proposition}[{\cite[Theorem 6.1, Corollary 6.1]{BKV}}]\label{prop:MainTool}
\begin{gather*}
\C^\bullet (C,M)\cong 
\C^\bullet ( \mathcal A_+(C), M)/\partial^*_\bullet \C^\bullet ( \mathcal A_+(C),M)\\
\tilde H^\bullet (C,M)= 
H^\bullet ( \mathcal A_+(C), M)
\end{gather*}
\end{proposition}

So we can calculate (reduced) Hochschild cohomologies of 
a conformal algebra $C$ via the Hochschild complex of its coefficient algebra $\mathcal A_+(C)$ and its quotient. 
To do that, we construct the Anick resolution for $\mathcal A_+(C)$ by means of the Morse matching method. 
We will use the homotopy maps constructed in this way 
to transfer the map $\partial^*_\bullet $
to the dual complex obtained from the Anick resolution,
as explained in the following sections.

\section{Morse matching and Anick resolution}\label{sec:Morse}

\subsection{The basics of algebraic discrete Morse theory}	
%
%
%
Suppose $\Bbbk $ is a field and $\Lambda $ is a unital associative  $\Bbbk$-algebra.
Let $\B_\bullet =(\B_n,\dcobound_n)_{n\ge 0}$ be a chain complex of free left $\Lambda $-modules, and let $X_n$ be a basis 
of $\B_n$ over $\Lambda $. One may represent the complex with a weighted oriented graph
$\Gamma (\B_\bullet )$
whose vertices are $\bigcup\limits_{n} X_n$ and there is an edge 
$\mathbf x \stackrel{\lambda} \longrightarrow \mathbf y $
from $\mathbf x\in X_n$ to $\mathbf y\in X_{n-1}$ 
of weight $\lambda $, $0\ne \lambda \in \Lambda $, 
if the distribution of $\dcobound_n(\mathbf x)\in \B_{n-1}$ as a linear 
combination of  $X_{n-1}$ contains $\lambda \mathbf y$.

A subset $M \subseteq E$ of the set of edges is called
Morse matching, if it satisfies the following three conditions:
\begin{itemize}
    \item  $M$ is a particular matching in the graph $\Gamma (\B_\bullet )$, i.e., a subset of edges such that 
neither of vertices belongs to more than one edge from $M$. 
\item The weights of edges from $M$ are invertible central elements of $\Lambda $. Then transform the graph $\Gamma(\B_\bullet )$ in the following way: invert the direction of all edges from $M$
and replace their weights $\lambda $ with $-\lambda^{-1}$.
\item Resulting graph $\Gamma_M(\B_\bullet )$ has no directed cycles.
\end{itemize}

 The vertices that do not belong to  edge from $M$ are said to be critical cells.

For a Morse matching graph $\Gamma_M(\B)$, one may construct another chain complex 
of free $\Lambda $-modules $(\A_m,\delta_m)_{m\ge 0}$, where $\A_m$ is spanned over 
$\Lambda $ by the set $X_{(m)}\subseteq X_m$ 
of all critical cells from $\B_m$ and the calculation of $\delta_m$ 
is based on the consideration of all paths in $\Gamma_M(\B)$  \cite{formancell,formanguide,JW}.

Namely,
\begin{equation}\label{eq:Anick}
\delta_m(\mathbf x)=\sum\limits_{\mathbf y\in  X_{(m-1)}}\Gamma(\mathbf x,\mathbf y)\mathbf y,\quad  \mathbf x\in  X_{(m)}
\end{equation}
where $\Gamma(\mathbf x,\mathbf y)$ is the sum of path weights in the Morse matching graph $\Gamma_M(\B)$.
The new complex $(\A_m,\delta_m)_{m\ge 0}$ is homotopy equivalent to the initial one, 
but it is smaller since we choose only critical cells as generators of the modules $\A_m$.

The homotopy maps $\mathrm f_n: \B_n \to \A_n$, $\mathrm g_n:\A_n \to \B_n$ are given by
\begin{align*}
     & \mathrm{f}_n(\mathbf {b}_n) = \sum_{\mathbf {a}_n\in \A_n} \Gamma_{\B_\bullet}(\mathbf {b}_n,\mathbf{a}_n) \mathbf{a}_n, \\
     & \mathrm{g}_n(\mathbf{a}_n) = \sum_{\mathbf{b}_n \in \B_n} \Gamma_{\A_\bullet}(\mathbf{a}_n,\mathbf{b}_n)\mathbf{b}_n.
 \end{align*}
If we are given a chain map $\varphi_\bullet : \B_\bullet \to \B_\bullet $
then for every $n\ge 1$ the map 
$\tilde \varphi _n = \mathrm f_n\varphi_n \mathrm  g_n$
is a chain map on $\A_\bullet$.
We have the following commutative diagram
\[
  \xymatrix{
   \A_n \ar@{->}[r]^{\delta_n} \ar@{->}[d]_{\mathrm{g}_n} & \A_{n-1} \\
   \B_n \ar@{->}[r]_{\dcobound_n} & \B_{n-1} \ar@{->}[u]_{\mathrm{f}_{n-1}}
  }
 \]
Therefore

\begin{equation}\label{com.}
\delta_n 
= \mathrm f_{n-1}\dcobound_n \mathrm g_n:
    \A_n\to \A_{n-1}.
\end{equation}

\subsection{The Anick resolution for associative algebras.}

Suppose $\Lambda $ has an an augmentation $\varepsilon: \Lambda \to \Bbbk $, 
$A$ is a set of generators for $\Lambda$. 
Then $\Lambda $ is a homomorphic image of the free associative algebra 
$\Bbbk \langle A\rangle$ generated by~$A$.
Assume  $\le$ is a monomial order on the free monoid $A^{*}$,
and $\GSB_\Lambda$ is a Gr\"obner---Shirshov basis of $\Lambda $.
The latter may be considered as a confluent set of defining relations 
for the algebra~$\Lambda $, each of relations is of the form $u-f$, where $u\in A^*$, 
$f\in \Bbbk \langle A\rangle$, $u\ge \bar f$, $\bar f$ is the leading monomial of $f$
relative to $\le $.
Denote by $V$ the set of all leading terms $u$ of relations from $\GSB_\Lambda$, $V$ is called the set of {\em obstructions}.

Following Anick \cite{Anick1983}, a word $v=x_{i_1}\ldots x_{i_t}$ 
is an $n$-{\em prechain} if and only if there exist $a_j,b_j \in \mathbb{Z}$, $1 \le j \le n$, satisfying the following conditions:
\begin{itemize}
    \item $1=a_1<a_2 \le b_1<\ldots<a_n \le b_{n-1}<b_n=t$;
    \item $ x_{i_{a_j}}\ldots x_{i_{b_j}} \in V$ for $1 \le j \le n$.
\end{itemize}
An $n$-prechain $x_{i_1}\ldots x_{i_t}$ is an {\em $n$-chain} 
if only if the integers $a_j,b_j$ can be chosen in such a way that  $x_{i_1}\ldots x_{i_t}$  is not an $m$-prechain for neither 
$s < b_m$, $1 \le m \le n$.

The set of all $n$-chains is denoted $V^{(n)}$.

The cokernel of  $\varepsilon: \Bbbk \to \Lambda$ is denoted by $\Lambda\slash \Bbbk$.
The set of all non-trivial words in $A^*$ that do not contain a word from $V$
as a subword (i.e., the set of $V$-reduced words \cite{Bok72})
forms a linear basis of $\Lambda/\Bbbk $. This is one of the equivalent 
conditions in the Composition-Diamond Lemma about Gr\"obner---Shirshov bases 
for associative algebras (see, e.g., \cite{Bokut}).
The resolution  $\B_\bullet =(\B_m,\dcobound_m)_{m\ge 0}$ is an 
exact sequence of $\Lambda$-modules where
$$
\B_m:=\Lambda \otimes(\Lambda/\Bbbk)^{\otimes_{m}}.
$$
and differential
$\dcobound_n : \B_n
   \to \B_{n-1}$.
   
We will use the standard convention denoting $1\otimes \alpha_1\otimes \dots \alpha_{m}\in \B_m $ 
 by $[\alpha_1|\ldots |\alpha_{m}]$ for $\alpha_i\in \Lambda/\Bbbk$ so differential is defined as follows:
\[
\dcobound_m([a_1| \ldots | a_n]) = a_1 
[a_2| \ldots |a_m]
+  \sum\limits_{i=1}^{m-1}(-1)^{i}   [a_1| \ldots |N(a_ia_{i+1})| \ldots | a_m].
\]
Here $N(a_ia_{i+1})$ is the corresponding Gr\"obner--Shirshov normal form of the product $a_ia_{i+1}$.

The {\em Anick resolution}  $\A_\bullet =(\A_m, \delta_m)_{m\ge 0}$ is an 
exact sequence of $\Lambda$-modules where
$$
\A_m:=\Lambda \otimes \Bbbk V^{(m-1)}.
$$
The role critical cells $X_{(m)}$ is played by  $(m-1)$-chain and it is not difficult to find it by $\GSB_\Lambda$. The computation of differentials in the Anick resolution  according to 
the original Anick algorithm described in \cite{Anick1983} is extremely hard.
In order to visualize the computation of differentials it is possible to use the discrete algebraic
Morse theory\cite{JW, formancell, formanguide}  based on the concept of a Morse matching.
For word $w$, let $\Lambda_{w,p}$ be set of all the vertices 
$[w_1| \ldots | w_n]$ in 
$\Gamma(\B_\bullet(\Lambda,\Bbbk))$ such that $w = w_1 \cdots w_n$ 
and~$p$ is the largest integer 
$p \ge -1$ for which $w_1\cdots w_{p+1} \in 
\Lambda^{(p)}$ is an Anick $p$-chain.
Let 
$\Lambda_w \coloneqq  \bigcup\limits_{p \ge -1}\Lambda_{w,p}$.

Define a partial matching $\mathcal{M}_w$ on 
$\Gamma(\B_\bullet(\Lambda,\Bbbk))|_{\Lambda_w}$ by letting 
$\mathcal{M}_w$ consist of all edges
\[
[w_1| \ldots| w'_{p+2}|w''_{p+2}| \ldots| w_n] \to [w_1| \ldots| w_{p+2}| 
\ldots| w_m]
\]
where 
$w'_{p+2}w''_{p+2} = w_{p+2}$, 
$[w_1|  \ldots| w_m] \in \Lambda_{w,p}$, 
and $[w_1| \ldots| w_{p+1}|w'_{p+2}] \in \Lambda^{(p+1)}$ is an Anick 
$(p+1)$-chain.

\section{The Anick complex for $\mathcal A_+(U(2))$}

\subsection{The Morse matching graph for $\mathcal A_+(U(2))$}\label{anick}
Let us write down the Gr\"obner--Shirshov basis of 
$\mathcal A_+ = \mathcal A_+(U(2))$ 
as of an associative algebra 
over a field $\Bbbk$ generated by
elements $v(n)$, $n\ge 0$.
The defining relations of $U(2)$ reflect 
the Virasoro commutator
$[v(n),v(m)]=(n-m)v(n+m-1)$ and
locality $N(v,v)=2$:
$v(n+2)v(m)-2v(n+1)v(m+1)+v(n)v(m+2)=0$, 
$n,m\ge 0$.

The Gr\"obner--Shirshov basis of $\mathcal A_+$
is easy to derive. It consists of the relations
\begin{equation}\label{eq:u3-gsb}
   v(n)v(m)=v(0)v(n+m)+nv(n+m-1),\quad n \ge 1, m \ge 0.
\end{equation}
Indeed, the linear basis of $\mathcal A_+ \simeq pW_1$
is given explicitly by the monomials $p^{k+1}q^n$, $k,n\ge 0$, 
that represent the reduced forms 
$v(0)^{k+1}v(n)$, $n\ge 0$, relative to \eqref{eq:u3-gsb}.

Throughout the rest of the paper 
$\Lambda $ stands for 
the augmented algebra 
$\Lambda = \mathcal A_+ \oplus \Bbbk $. 

%

To find the Anick complex, we need two steps.
First, we have to find the set of obstructions 
for $\mathcal A_+ $ 
relative to the given Gr\"obner---Shirshov basis (the set of leading terms in $\mathcal A_+(U(2))$ ) 
and the set of $n$-chains. 
Next, build a Morse graph and calculate the path weights 
$\Gamma(w,w')$ for every $n$-chain $w$ and $(n-1)$-chain $w'$ for all $n\ge1$. For all $n\ge2$ we have
\[
V^{(n)}=\{v(i_1)v(i_2)\ldots v(i_{n+1});\quad i_1,i_2,\ldots,i_n\ge1,i_{n+1}\ge0\}.
\]

For $n=1$ we have the following set of obstructions
\[
V^{(1)}=\{v(n)v(m);\quad n\ge1,m\ge0\}.
\]
Let us  
evaluate $\delta_2: \A_{2}\to \A_{1}$ by means of the Morse graph is shown on Figure \ref{FIG:1}.

Hence,
\begin{figure} 
\centering
\begin{tikzpicture}[commutative diagrams/every diagram]
\node (A1) at (0,0) {$[v(n)|v(m)]$};

\node (B1) at (-4,-1) {$[v(m)]$};
\node (B2) at (-2, -2) {$[v(0)v(n+m)]$};
\node (B3) at (3,-1.5)  {$[v(n+m-1)]$};

\node (C1) at (-2,-4) {$[v(0)|v(n+m)]$} ;

\node (D1) at (-4, -5.5) {$[v(n+m)]$};

\path[commutative diagrams/.cd,every arrow, every label]
  (A1) edge node [swap] {$v(n)$} (B1)
  (A1) edge node [swap] {$-1$} (B2)
  (A1) edge node [right] {$-n$} (B3)
  (C1) edge node [swap] {$v(0)$} (D1)
   ;

 \path[->]
  ([xshift=5pt]C1.north) edge node[right] {$-1$} ([xshift=5pt]B2.south);
  \path[->,dashed]
  ([xshift=-5pt]B2.south) edge node[left]{$1$} ([xshift=-5pt]C1.north);

\end{tikzpicture}
    \caption{}
    \label{FIG:1}
\end{figure}
\[
\delta_2[v(n)|v(m)]=v(n)[v(m)]-v(0)[v(n+m)]-n[v(n+m-1)]
\]
For $n=2$ we have
\[
V^{(2)}=\{v(n)v(m)v(p);\quad n,m\ge1,p\ge0\}.
\]
and $\delta_3: \A_{3}\to \A_{2}$ is given by the following rule:
\[
\begin{aligned}
  \delta_3[v(n)v(m)v(p)]{} = v(n)[v(m)v(p)]- n[v(n+m-1)v(p)] - v(0)[v(n+m)v(p)] \\+ v(0)[v(n)v(m+p)] + n[v(n-1)v(m+p)]
			+ m[v(n)v(m+p-1)];\quad n,m\ge1,p\ge0.
\end{aligned}
\]

\begin{figure} 
\centering
 \begin{tikzpicture}[commutative diagrams/every diagram]
  \node (A) at (0,0) { $[v(i_1)|\cdots | v(i_{n+1})]$ };
  \node (B0) at (-5,-1) {$[v(i_2)|\cdots | v(i_{n+1})]$};
  \node (C1) at (-4,-5) {$[ v(i_1)| \cdots | v(0)v(i_j+i_{j+1})| \cdots v(i_{n+1}) ]$};
  \node (C2) at (4,-5) {$[v(i_1)|\cdots | v(i_j + i_{j+1}-1)| \cdots v(i_{n+1})]$};
  
  \path[commutative diagrams/.cd,every arrow, every label]
  (A) edge node [swap] {$v(i_1)$} (B1)
  (A1) edge node [swap] {$(-1)^j$} (C1)
  (A1) edge node {$(-1)^ji_j$} (C2)
  ;
\end{tikzpicture}
    \caption{}
    \label{fig:3}
\end{figure}
\begin{figure} 
\centering
\begin{tikzpicture}[commutative diagrams/every diagram]
\node (A) at (0,0) {$[v(i_1)|\cdots | v(i_{t-1})|v(0)v(p)|\cdots| v(i_{n+1})]$};
\node (B) at (0,-2) {$[v(i_1)|\cdots| v(i_{t-1})|v(0)|v(p)| \cdots| v(i_{n+1})]$};
\node(B0) at(-5,-2) {$\substack{ \mbox{no} \\ \mbox{Anick} \\ \mbox{chains} }$};
\path[->]
  ([xshift=-10pt]A.south) edge node[left] {$(-1)^{t+1}$} ([xshift=-10pt]B.north);
\path[->,dashed]
 ([xshift=10pt]B.north) edge node[right]{$(-1)^{t}$} ([xshift=10pt]A.south);
\path[->]
 ([yshift=5pt]B.west) edge node[yshift=-4pt] {$\vdots$} ([yshift=9pt]B0.east);
\path[->]
 ([yshift=-5pt]B.west) edge node {} ([yshift=-9pt]B0.east);
\end{tikzpicture}
    \caption{}
    \label{fig:4}
\end{figure}
we can evaluate Anick differential $\delta_{n}: \A_{n}\to \A_{n-1}$ as we found $\delta_3$ by repeating the steps \cite{A2022}, which are shown in the two Figure \ref{fig:3}, \ref{fig:4} 
then $\delta_n$ is given by the following rules:
\[
\begin{aligned}
	&\delta_{n}[v(i_1)v(i_2)|\ldots|v(i_{n-1})v(i_n)]= v(i_1)[v(i_2)v(i_3)\ldots v(i_{n-1})v(i_n)]\\
	&+\sum\limits_{j=1}^{n-1}(-1)^ji_j[v(i_1)|v(i_2)|\ldots|v(i_j+i_{j+1}-1)|\ldots|v(i_n)] +\\ &+\sum\limits_{j=1}^{n-1}(-1)^{j}v(0)[v(i_1)|v(i_2)|\ldots|v(i_j+i_{j+1})|\ldots|v(i_n)] + \\
	&+\sum\limits_{j=2}^{n-1}\sum\limits_{k=1}^{j-1}(-1)^{j}i_k[v(i_1)|\ldots|v(i_k-1)|\ldots|v(i_j+i_{j+1})|\ldots|v(i_n)];\quad i_1,\ldots,i_{n-1}\ge1,i_n\ge0.
	\end{aligned}
\]

\subsection{Anick resolution for differential algebras}

Let $C$ be an associative conformal algebra. Throughout the rest of the paper 
$\Lambda $ stands for 
the augmented algebra  
$\Lambda = \mathcal A_+ \oplus \Bbbk 1$,
where 
$\mathcal A_+ = \mathcal A_+(C)$, 
and the augmentation is given by
$\varepsilon (\mathcal A_+)=0$.
Then $\mathcal A_+=\Lambda /\Bbbk $, 
and
$\Lambda $ has a derivation $\partial $ such that
$\partial (a(n)) = -na(n-1)$, $\partial (1)=0$. 


For every $n\ge 1$, there is a $\Bbbk$-linear map 
\[
\partial_n: \B_n\to \B_n,
\]
\[
\partial_n(\beta [\alpha_1|\dots |\alpha_n])
= \partial(\beta) [\alpha_1|\dots |\alpha_n]
+\sum\limits_{i=1}^n 
 \beta [\alpha_1 | \dots |\partial(\alpha_i)| \dots |\alpha_n],
\]
$\alpha_i\in \mathcal A_+$, $\beta \in \Lambda $. 
This is not a morphism of complexes, 
but it induces a morphism of dual cochain complexes
which can be transferred to the Anick resolution 
via homotopy.
Namely, suppose $M$ is a left conformal $C$-module. Thus, 
$M$ is an ``ordinary'' $\Lambda $-module. 
Denote
\[
\mathrm C_{\B}^n = \Hom _\Lambda (\B_n, M) \simeq 
\Hom_\Bbbk (\A^{\otimes n}, M), 
\]
this is the space of Hochschild cochains. 
The Hochschild differential 
\[
\Delta_{\B}^n  : \C_{\B}^n \to \C_{\B}^{n+1} 
\]
is expressed as 
\[
\big (\Delta_{\B}^n\varphi\big )(\mathbf x) = \varphi \dcobound_{n+1}(\mathbf x), \quad \varphi\in \C_{\B}^{n},
\  \mathbf x\in \B_{n+1}.
\]

Note that the $\Lambda $-module $M$ is equipped with a derivation 
also denoted $\partial $ (the same as in the definition of a conformal module), so that 
$\partial (a(n)u) = -na(n-1)u +a(n)\partial(u)$, for $a\in \C$, $u\in M$, 
$n\in \mathbb Z_+$.

Then for every $n\ge 1$ the map
\[
D_{\B}^n : \C_{\B}^n \to \C_{\B}^n 
\]
given by
\[
(D_{\B}^n \varphi)(\mathbf x) = \partial(\varphi (\mathbf x))
-\varphi (\partial_n(\mathbf x)),
\quad 
\varphi\in \mathrm C_{\B}^{n},
\  \mathbf x\in \B_{n+1},
\]
is a morphism of complexes:
\[
D_{\B}^{n+1}\Delta^n_{\B} 
=  \Delta^n_{\B} D_{\B}^{n}.
\]

Let us translate the mapping $D^\bullet _{\B}$ 
from 
$\C_{\B}^\bullet $
to the complex 
$\C_{\A}^\bullet $
constructed on the spaces
\[
\C_{\A}^n = \Hom_\Lambda (\A_n, M)\simeq \Hom_{\Bbbk }
(\Bbbk V^{(n-1)}, M)
\]
with the differential $\Delta_{\A}^n:  \C_{\A}^n \to \C_{\A}^{n+1}$ given by (see \ref{com.})
\[
\Delta^n_{\A}(\psi) = \psi\delta_{n+1}
= \psi \mathrm f_n\dcobound_{n+1} \mathrm g_{n+1}
=\Delta_{\B}^n(\psi \mathrm f_n) g_{n+1}.
\]
The homotopy equivalence between $\A_\bullet $
and $\B_\bullet $
turns $D_{\B}^n$
into 
\[
D_{\A}^n : \C_{\A}^n \to \C_{\A}^{n}
\]
such that 
\[
D_{\A}^n\psi = D_{\B}^n(\psi \mathrm f_n)\mathrm g_n.
\]
For every $\mathbf a\in \A_n$, 
we have 
\[
(D_{\mathrm A}^n\psi)(\mathbf a)
= \partial (\psi \mathrm f_n (\mathrm g_n (\mathbf a)))
-(\psi\mathrm f_n)(\partial_n(\mathrm g_n(\mathbf a)))
= \partial (\psi (\mathbf a))
-(\psi\mathrm f_n)(\partial_n(\mathrm g_n(\mathbf a))).
\]

Note that the Anick chains for $\mathcal A_+=\mathcal A_+(U(2))$
constructed in Section~\ref{anick}
have the following property: if $\mathbf x = [\alpha_1|\dots |\alpha_n]\in 
\mathrm B_n$
is not an Anick chain then $\partial_n(\mathbf x)$ does not contain 
summands that are Anick chains, i.e.,
$\mathrm f_n(\mathbf x)=0$ implies
$\mathrm f_n(\partial_n(\mathbf x)) = 0$.
Hence, we may evaluate $D_{\mathrm A}^n\psi $ on $\mathbf a\in \mathrm A_n $
in the same way as 
$D_{\mathrm B}^n(\psi)$ on $\mathbf a$, 
simply removing those terms from $\mathrm B_n$ that do not represent 
Anick chains.

\begin{example}
Let $\psi \in \mathrm C_{\mathrm A}^3$ for $A=\mathcal A_+(U(2))$.
Denote
$\psi([v(n_1)|v(n_2)|v(m)]) = \psi(n_1,n_2,m)\in M$,
for $n_1,n_2\ge 1$, $m\ge 0$.
Let us calculate 
$(D^3_{\mathrm A}\psi )[v(2)|v(1)|v(1)]$.
\end{example}

According to the Morse matching graph in Sec \ref{sec:Morse},
we have 
\begin{multline*}
\mathrm g_3[v(2)|v(1)|v(1)] = [v(2)|v(1)|v(1)]-[v(0)|v(3)|v(1)]-[v(2)|v(0)|v(2)] \\
+[v(0)|v(0)|v(4)]+[v(2)|v(0)|v(2)].
\end{multline*}
Then 
\[
\partial_3(g_3[v(2)|v(1)|v(1)]) = -2[v(1)|v(1)|v(1)]-[v(2)|v(1)|v(0)] + \dots,
\]
where all other summands are not Anick chains (contain $v(0)$ not at the last position).
Therefore, 
\[
(D^3_{\mathrm A}\psi )[v(2)|v(1)|v(1)] = \partial \psi(2,1,1) + 2\psi(1,1,1)+\psi(2,1,0).
\]

Proposition \ref{prop:MainTool} and homotopy equivalence of 
$\B_\bullet$ and $\A_\bullet$
immediately imply the following statement.

\begin{proposition}
  For a conformal algebra $C$ and a conformal $\C-$module $M$, 
  the reduced Hochschild complex 
  $\C^\bullet (\C,M)$ is homotopy equivalent 
  to $\C_{\A}^\bullet / D_{\A}^\bullet \C_{\A}^\bullet $.
\end{proposition}

The complex
$\mathrm C_{\mathrm A}^\bullet /D_{\mathrm A}^\bullet \mathrm C_{\mathrm A}^\bullet $ constructed 
on the algebra $\mathcal A_+ = \mathcal A_+(C)$
is called the {\em reduced Anick complex} of~$\mathcal A_+$.

\section{Cohomologies of the reduced Anick complex of $\mathcal A_+(U(2))$}

Given a conformal module $M$ over the 
associative conformal algebra $U(2)$, we denote 
by $\C^\bullet $
 the reduced Anick complex for 
$\mathcal A_+ = \mathcal A_+(U(2))$, and let
 $\tilde \C^\bullet$ stand for the non-reduced 
complex $\C_{\mathrm A}^\bullet $ with values in~$M$.
By $D^n$ we will denote the operation $D^n_{\mathrm A}$
on $\tilde C^n$.
In order to simplify notations, we will write 
$[i_1|i_2|\dots |i_n]$ for $[v(i_1)|v(i_2)| \ldots | v(i_n)]\in \mathrm B_n$.

In the case of non-scalar module $M$, like $M_{(\alpha,\Delta)}$, 
it is technically easier to work with the reduced complex 
$\C^\bullet $ rather than 
with the non-reduced one $\tilde{\C}^\bullet $.

\begin{lemma}\label{lem:Reduction}
The elements of $\C^n$ are in one-to-one 
correspondence with 
scalar sequences $\alpha_{(i_1,i_2,\ldots,i_n)}$, 
$[i_1|i_2|\ldots |i_n]\in \mathrm A_n$.
\end{lemma}

\begin{proof}
Let us identify elements of type $f(\partial )u \in M_{(\alpha,\Delta )}$ with polynomials $f(\partial )\in \Bbbk [\partial ]$.

Introduce the lexicographic order $\le_{\text{lex}}$ 
on the Anick chains 
$[i_1|i_2|\ldots |i_{n-1}|i_n]\in \mathrm A_n$ 
as on the sequences of non-negative integers.
This is a well order with the smallest element
$[1|1|\ldots |1|0]$. Let us prove by induction 
that 
for every $f\in \tilde \C^n$
there exists unique sequence of polynomials 
\[
h_{(j_1,j_2,\ldots ,j_{n-1},j_n)} \in \Bbbk [\partial ], 
\quad 
[j_1|j_2|\ldots |j_{n-1}|j_n] \in \mathrm A_n,
\]
such that 
\begin{equation}\label{eq:DerivationCoChain}
f[i_1|i_2|\ldots |i_{n-1}|i_n] - 
\partial h_{(i_1,i_2,\ldots ,i_{n-1},i_n)}
-\sum\limits_{k=1}^n i_k h_{(i_1,i_2,\ldots , i_k-1 , \ldots,i_{n-1},i_n)} = \alpha_{(i_1,i_2,\ldots,i_{n-1},i_n)}
\in \Bbbk 
\end{equation}
for every 
$[i_1|i_2|\ldots |i_{n-1}|i_n]\in \mathrm A_n$.
Here we assume that 
$h_{(j_1,j_2,\ldots ,j_{n-1},j_n)}=0$
if 
$[j_1|j_2|\ldots |j_{n-1}|j_n] \notin \mathrm A_n$.
If \eqref{eq:DerivationCoChain} holds then 
$f-D^nh$ takes scalar values at $\mathrm A_n$
for 
$h\in \tilde \C^n$ given by 
\[
h[i_1|i_2|\ldots |i_{n-1}|i_n] = h_{(i_1,i_2,\ldots ,i_{n-1},i_n)}.
\]

It is easy to find $h_{(1,1,\ldots, 1,0)}$ as 
$(b(\partial )- b(0))/\partial $
for $b = f[1|1|\ldots |1|0]$.

Assume $h_{(j_1,j_2,\ldots ,j_{n-1},j_n)}$
are already constructed for all 
$[j_1|j_2|\ldots |j_{n-1}|j_n]\le_{\text{lex}}
[i_1|i_2|\ldots |i_{n-1}|i_n]$.
Then 
$h_{(i_1,i_2,\ldots ,i_{n-1},i_n)}$
is uniquely defined 
as  $(b(\partial )- b(0))/\partial $
for 
\[
b(\partial )=
f[i_1|i_2|\ldots |i_{n-1}|i_n] 
+\sum\limits_{k=1}^n i_k h_{(i_1,i_2,\ldots , i_k-1 , \ldots,i_{n-1},i_n)} .
\]
\end{proof}

\begin{theorem} For the conformal module $M_{(\alpha,1)}$ 
over $U(2)$, we have 
\[ 
\dim_\Bbbk\Homol^1(U(2), M_{(\alpha,1)})=
\begin{cases}
1 & \alpha=0 ,\\
0 & \alpha\ne0.
\end{cases}
\]
\end{theorem}

\begin{proof}
We are interested in the space
$\Homol^1(U(2),M) $ which is isomorphic 
to the space of non-coboundary cocycles in $\C^1 = \tilde \C^1/D^1\tilde \C^1$.
Suppose $\varphi \in \tilde \C^1 = \Hom_\Lambda (\A_1, M)$. 
By Lemma~\ref{lem:Reduction} we may assume 
$\varphi [n] = \alpha_n \in \Bbbk $
for $n\ge 0$.
Then the differential $\Delta^1\varphi$ takes the following 
values on the Anick 1-chains (see Figure~\ref{FIG:1}):
\[
(\Delta^1\varphi)[n|m] =\varphi(\delta_2[n|m]) = v(n)\alpha_{m}-v(0)\alpha_{n+m}-n\alpha_{n+m-1}.
\]
Hence by \eqref{eq:v(n)-Action} we have 
\[
\begin{gathered}
(\Delta^1\varphi)[1|m] =-(\partial+\alpha)\alpha_{1+m},\quad m\ge0, 
\\
(\Delta^1\varphi)[n|m]  =-(\partial+\alpha)\alpha_{n+m}-n\alpha_{n+m-1},\quad n\ge2,\ m\ge0.
\end{gathered}
\]
In order to find the constants that define $\Delta^1(\varphi +D^1\tilde \C^1) \in \C^2$
choose $\psi \in \tilde C^2$
such that 
\[
 \psi[n|m] =- \alpha_{n+m},\quad n\ge1,\ m\ge0.
\]
Then  
\[
(\Delta^1\varphi-D^2\psi)[n|m]=-\alpha\alpha_{n+m}+m\alpha_{n+m-1}
\]
for all $[n|m]\in \mathrm A_2$.
Therefore, $\varphi +D^1\tilde \C^1 $ is a 1-cocycle in $\C^1$
if and only if 
$\alpha \alpha_{n+m} = m \alpha_{n+m-1}$
for all $n\ge 1$, $m\ge 0$.
The latter is possible only in $\alpha_m=0$ for all $m\ge 1$.
Hence, 1-cocycles in $\C^1$ form a 1-dimensional space:
every 1-cocycle is determined by $\alpha _0$.

On the other hand, coboundary cocycles in $\C^1$ are given by 
$\Delta^0 (h+D^0\tilde C^0)$, where $h\in \Bbbk [\partial ]$. 
Modulo $D^0\tilde C^0$, we may assume $h(0) = \beta \in \Bbbk $, 
then 
\[
(\Delta^0h)[n] = v(n)\beta =
\begin{cases}
(\partial + \alpha)\beta , & n=0, \\
\beta , & n=1, \\
0, & n>1. \end{cases}
\]
Reduce $\partial $ modulo $D^1\tilde \C^1$: 
choose $\xi \in \tilde \C^1$ such that 
$\xi [n]=\delta_{n,0}\beta  $.
Then $(D^1\xi )[n] = \delta_{n,0}\partial \beta + \delta_{n,1}\beta $,
\[
(\Delta^0 h - D^1\xi )[n]=\delta_{n,0}\alpha\beta.
\]
As a result, the space of 1-coboundaries is 1-dimensional for 
$\alpha\ne 0$ and zero for $\alpha =0$.
Hence, 
\[ 
\dim_\Bbbk\Homol^1(U(2), M_{(\alpha,1)})=
\begin{cases}
1 & \alpha=0,\\
0 & \alpha\ne0.
\end{cases}
\]
\end{proof}

The following statement is a corollary of 
\cite{Kozlov2017} where it was proved that 
the second Hochschild cohomology group of $U(2)$
is always trivial. Let us still present its proof 
by means of our methods since the same approach 
would work later for an arbitrary $n\ge 2$.

\begin{theorem}\label{thm:Anick} The the second
cohomology group of $U(2)$ with the value in $M_{(\alpha,1)}$ is trivial, 
\[
\Homol^2(U(2),M_{(\alpha,1)})=0.
\]
\end{theorem}

\begin{proof}
 Suppose $\varphi \in \tilde \C^2$. 
The class $\varphi +D^2\tilde \C^2$ is completely 
defined by a sequence of scalars $\alpha_{(n,m)}$
by Lemma~\ref{lem:Reduction}.

The differential $\Delta^2\varphi$ takes the following 
values on the Anick 2-chains ($n,m\ge 1$, $p\ge 0$):
\begin{multline*}
(\Delta^2\varphi)[n|m|p] = \varphi (\delta_3[n|m|p])
= v(n)\alpha_{(m,p)} - n\alpha_{(n+m-1,p)}- v(0)\alpha_{(n+m,p)}\\
   + v(0)\alpha_{(n,m+p)} + n\alpha_{(n-1,m+p)}	+ m\alpha_{(n,m+p-1)}.
\end{multline*}
The summand $n\alpha_{(n-1,m+p)}$ does not appear if $n=1$ since 
$[0|m+p]$ is not an Anick chain. Hence, we have 
to consider two cases: $n=1$ and $n\ge 2$.

For all $n\ge2$, $m\ge1$, $p\ge0$, we have 
\[
    (\Delta^2\varphi)[n|m|p] 
    = - n\alpha_{(n+m-1,p)}- (\partial+\alpha)\alpha_{(n+m,p)}
   + (\partial+\alpha)\alpha_{(n,m+p)}
   + n\alpha_{(n-1,m+p)}	+ m\alpha_{(n,m+p-1)}
\]
by \eqref{eq:v(n)-Action}.
Similarly, for $n=1$, $m\ge1$, $p\ge0$ we obtain
\[
(\Delta^2\varphi)[1|m|p] =- (\partial+\alpha)\alpha_{(1+m,p)} + (\partial+\alpha)\alpha_{(1,m+p)}+ m\alpha_{(1,m+p-1)}. 
\]
Reduce the result by means of $D^3\psi $, where
$\psi \in \tilde C^3 $ is given by 
\[
\psi[n|m|p]=-\alpha_{(n+m,p)} + \alpha_{(n,m+p)}, 
 \quad n\ge1,\ m\ge1,\ p\ge0.
    \]
In both cases, we obtain the same expression for 
 $\Delta^2\varphi - D^3\psi$:
\[
  (\Delta^2\varphi-D^3\psi)[n|m|p] =
  - \alpha\alpha_{(n+m,p)} + \alpha\alpha_{(n,m+p)} 
  + m\alpha_{(n+m-1,p)} + p\alpha_{(n+m,p-1)}- p\alpha_{(n,m+p-1)}.
\]
Hence, $\varphi +D^2\tilde \C^2$ is a 2-cocycle in $\C^2$
if and only if 
\begin{equation}\label{eq:2-cocycleCond}
 - \alpha\alpha_{(n+m,p)} + \alpha\alpha_{(n,m+p)} 
  + m\alpha_{(n+m-1,p)} + p\alpha_{(n+m,p-1)}- p\alpha_{(n,m+p-1)} = 0
\end{equation}
for all $[n|m|p]\in \mathrm A_3$.

{\sc Case 1:} $\alpha \ne 0$. 
Put $p=0$ in \eqref{eq:2-cocycleCond} to obtain  
\[
 \begin{gathered}
-\alpha\alpha_{(1,m)}= - \alpha\alpha_{(1+m,0)}+m\alpha_{(m,0)}, \quad m\ge1, \\
-\alpha\alpha_{(n,m)}= - \alpha\alpha_{(n+m,0)}+m\alpha_{(n+m-1,0)}, 
\quad n\ge2,\ m\ge1.
  \end{gathered}
\]
Therefore, cocycles in $\C^2$ are determined by $\alpha_{(n,0)}$
for $n\ge 1$. 

Choose $\varphi_1\in \tilde \C^1$, $\psi_1 \in \tilde \C^2$ 
such that
    \[
    \begin{aligned}
    \varphi_1[n]&= 
    -\dfrac{\alpha_{(n,0)}}{\alpha},\quad n \ge 1.\\
    \psi_1[n|m]& = - \alpha_{n+m}, \quad n\ge1,\  m\ge0.
    \end{aligned}
    \]
Then $\Delta^1\varphi_1 $ is a coboundary in $\tilde \C^2$, and
\[
\begin{gathered}
(\Delta^1\varphi_1-D^2\psi_1)[n|0]=\alpha_{(n,0)}=\varphi[n|0], \quad n\ge1.
   \end{gathered}
\]
Hence, 
$\varphi -\Delta^1\varphi_1 \in D^2\tilde \C^2$, so 
every cocycle in $\C^2$ is a coboundary.

{\sc Case 2:} $\alpha =0$. 
 In this case, we have
\begin{equation}\label{eq:2-cocycleCond-II}
 m\alpha_{(n+m-1,p)} + p\alpha_{(n+m,p-1)}- p\alpha_{(n,m+p-1)} = 0
\end{equation}
for all $[n|m|p]\in \mathrm A_3$.

Put $p=1$ in \eqref{eq:2-cocycleCond-II} to get
   \[
   0=m\alpha_{(n+m-1,1)}+\alpha_{(n+m,0)}- \alpha_{(n,m)},
   \quad n\ge 1,\ m\ge1.
   \]
For $n=1$ and $p=0$ we simply obtain from \eqref{eq:2-cocycleCond-II}
that $\alpha_{(m,0)}=0$  for all $m\ge 1$.
So 
   \[
    \alpha_{(n,m)}=m\alpha_{(n+m-1,1)},\quad n\ge1,\ m\ge1,
    \]
i.e., cocycles in $C^2$ are determined by 
 $\alpha_{(n,1)}$, $n\ge1$.

Choose 
$\varphi_1\in\tilde \C^1$ and $\psi_1\in\tilde \C^2$ such that
    \[
    \begin{aligned}
    \varphi_1[n]&=\alpha_{(n,1)},\quad n\ge1, \\
  \psi_1[n|m]& =- \alpha_{n+m},\quad n\ge1, \ m\ge0.
    \end{aligned}
    \]
Then $\Delta^1\varphi_1 $ is a coboundary in $\tilde \C^2$, 
and
\[
(\Delta^1\varphi_1-D^2\psi_1)[n|1]=\alpha_{(n,1)},
\quad n\ge1.
\]
Hence, every 2-cocycle is a coboundary.
\end{proof}

\begin{theorem}\label{the5.5}
 For a conformal module $M_{(\alpha,1)}$ over the 
associative conformal algebra $U(2)$, and for all $n\ge3$  we have
 \[
 \Homol^{n-1}(U(2),M_{(\alpha,1)})=0 .
 \]
\end{theorem}
 
 \begin{proof}
Suppose $\varphi \in \tilde \C^{n-1} = \Hom_\Lambda (\A_{n-1}, M)$, $M=M_{(\alpha, 1)}$. 
As above, let us identify $f(\partial )u\in M$
and $f\in \Bbbk [\partial ]$.
Suppose
$\varphi[i_1|i_2|\ldots |i_{n-1}] = f_{(i_1,i_2,\ldots,i_{n-1})}(\partial) \in \Bbbk [\partial]$ for $[i_1|i_2|\ldots |i_{n-1}]\in \mathrm A_{n-1}$.
It is convenient to assume that 
$f_{(i_1,i_2,\ldots,i_{n-1})}(\partial)=0$
if 
$[i_1|i_2|\ldots |i_{n-1}]\notin \mathrm A_{n-1}$.

Recall that the operation $D^{n-1}$
acts on $\tilde \C^{n-1}$
as 
\[
(D^{n-1}\varphi )[i_1|i_2|\ldots|i_{n-1}] = \partial f_{(i_1,i_2,\ldots,i_{n-1})} + i_1f_{(i_1-1,i_2,\ldots,i_{n-1})}+\ldots+i_{n-1}f_{(i_1,i_2,\ldots,i_{n-1}-1)}. 
\]
By Lemma~\ref{lem:Reduction}, the elements of $\C^{n-1}=\tilde \C^{n-1}/D^{n-1}\tilde \C^{n-1}$ 
are defined by 
scalar sequences 
$\alpha_{(i_1,i_2,\ldots,i_{n-1})}$.

The differential $\Delta^{n-1}$ takes the following
values on the Anick $(n-1)$-chains:
\[
\begin{aligned}
	&(\Delta^{n-1}\varphi)[1|i_2|\ldots|i_{n-1}|i_n]=\alpha_{(i_2,i_3,\ldots,i_{n-1})}\\
	&+\sum\limits_{j=1}^{n-1}(-1)^ji_j\alpha_{(1,i_2,\ldots,i_j+i_{j+1}-1,\ldots,i_{n})} +\\ &+\sum\limits_{j=1}^{n-1}(-1)^j(\partial+\alpha)\alpha_{(1,i_2,\ldots,i_j+i_{j+1},\ldots,i_n)} + \\
	&+\sum\limits_{j=2}^{n-1}\sum\limits_{k=1}^{j-1}(-1)^ji_k\alpha_{(1,\ldots,i_k-1,\ldots,i_j+i_{j+1},\ldots,i_n)},
 \quad i_2,\ldots,i_{n-1}\ge1, \ i_n\ge0,
	\end{aligned}
\]
\[
\begin{aligned}
	&(\Delta^{n-1}\varphi)[i_1|i_2|\ldots|i_{n-1}|i_n]=
	\sum\limits_{j=1}^{n-1}(-1)^ji_j\alpha_{(i_1,i_2,\ldots,i_j+i_{j+1}-1,\ldots,i_{n})} \\ &+\sum\limits_{j=1}^{n-1}(-1)^j(\partial+\alpha)\alpha_{(i_1,i_2,\ldots,i_j+i_{j+1},\ldots,i_n)}\\
	&+\sum\limits_{j=2}^{n-1}\sum\limits_{k=1}^{j-1}(-1)^ji_k\alpha_{(i_1,\ldots,i_k-1,\ldots,i_j+i_{j+1},\ldots,i_n)},
 \quad i_1\ge2,\ i_2,\ldots,i_{n-1}\ge1,\ i_n\ge0.
	\end{aligned}
\]
 Let us reduce 
 $\Delta^{n-1}\varphi $
 by means of $D^n\psi $ for 
 $\psi \in \tilde C^{n} $
 given by
\[
\psi[i_1|i_2|\ldots|i_{n}]
=\sum\limits_{j=1}^{n-1}(-1)^j\alpha_{(i_1,i_2,\ldots,i_j+i_{j+1},\ldots,i_n)}.
\]
Namely, 
\begin{align}
\nonumber
  (\Delta^{n-1}\varphi-D^n\psi)[i_1|i_2|\ldots|i_{n}]&
  {} =\sum\limits_{j=1}^{n-1}(-1)^{j}\alpha\alpha_{(i_1,i_2,\ldots,i_j+i_{j+1},\ldots,i_n)}\\
 &
\nonumber
 +\sum_{j=3}^{n}\sum_{k=1}^{j-2}(-1)^{k}i_j\alpha_{(i_1,\ldots,i_k+i_{k+1},\ldots,i_j-1,\ldots,i_n)}\\
 &
+\sum^{n-1}_{j=1}(-1)^{j+1}i_{j+1}\alpha_{(i_1,\ldots, i_j+i_{j+1}-1,\ldots ,i_n)} .
\label{eq:Cocycle-n-1}
\end{align}
Suppose $\varphi + D^{n-1}\tilde \C^{n-1} $
is a cocycle in $\C^{n-1}$. 
Then the right-hand side of 
\eqref{eq:Cocycle-n-1} is zero for every 
$[i_1|i_2|\ldots |i_n]\in \mathrm A_n$.

{\sc Case 1:} $\alpha\ne 0$.
Then 
for $i_1,i_2,\ldots,i_{n-1}\ge1$, $i_n=0$
relation \eqref{eq:Cocycle-n-1} implies
\begin{multline}\label{eq:5.4}
 (-1)^{n}\alpha\alpha_{(i_1,i_2,,\ldots,i_{n-1})} =\sum\limits_{j=1}^{n-2}(-1)^{j}\alpha\alpha_{(i_1,i_2,\ldots,i_j+i_{j+1},\ldots,i_{n-1},0)}\\
 +\sum_{j=3}^{n-1}\sum_{k=1}^{j-1}(-1)^{k+1}i_j\alpha_{(i_1,\ldots,i_k+i_{k+1},\ldots,i_j-1,\ldots,i_{n-1},0)}\\
  +\sum^{n-2}_{j=1}(-1)^{j+1}i_{j+1}\alpha_{(i_1,\ldots,
  i_j+i_{j+1}-1,\ldots ,i_{n-1},0)}
\end{multline}
Hence, the entire sequence 
$\alpha_{(i_1,i_2,\ldots,i_{n-1})}$
is determined by $\alpha_{(i_1,i_2,,\ldots,i_{n-2},0)}$ for $i_1,i_2,\ldots,i_{n-1}\ge 1$.

Choose $\varphi_1\in \tilde C^{n-2}$, 
$\psi_1\in \tilde \C^{n-1}$ 
such that
\[
\varphi_1[i_1|i_2|\ldots |i_{n-2}]
 =\beta_{(i_1,i_2,\ldots,i_{n-2})},
 \quad 
 \psi_1[i_1|i_2|\ldots|i_{n-1}] 
 =\sum\limits_{j=1}^{n-2}(-1)^{j}
  \beta_{(i_1,i_2,\ldots,i_j+i_{j+1},\ldots,i_{n-1})},
\]
where 
\[
\beta_{(i_1,i_2,\ldots,i_{n-2})}  
=
\begin{cases}
(-1)^{n}\dfrac{\alpha_{(i_1,i_2,\ldots,i_{n-2},0)}}{\alpha}, & i_1,i_2,\ldots,i_{n-1}\ge1, \\
0,& \text{otherwise}.
\end{cases}
\]
Then $\Delta^{n-2}\varphi_1-D^{n-1}\psi_1$
represents a coboundary in $C^{n-1}$, and
\[
\begin{aligned}
  (\Delta^{n-2}\varphi_1-D^{n-1}\psi_1)
   [i_1|i_2|\ldots |i_{n-1}] 
   & 
   =\sum\limits_{j=1}^{n-2}(-1)^{j}\alpha
      \beta_{(i_1,i_2,\ldots,i_j+i_{j+1},\ldots,i_{n-1})}
      \\
 &+\sum_{j=3}^{n-1}\sum_{k=1}^{j-2}(-1)^{k}i_j
  \beta_{(i_1,\ldots,i_k+i_{k+1},\ldots,i_j-1,\ldots,i_{n-1})}
  \\
     &+\sum^{n-2}_{j=1}(-1)^{j+1}i_{j+1}
 \beta_{(i_1,\ldots ,i_j+i_{j+1}-1,\ldots ,i_{n-1})}.
\end{aligned}
\]
For all $i_1,i_2,\ldots,i_{n-1}\ge1$, we get
\[
(\Delta^{n-2}\varphi_1-D^{n-1}\psi_1)
[i_1|i_2|\ldots |i_{n-2}|0]
 =\alpha_{(i_1,i_2,\ldots,i_{n-2},0)}.
\]
Hence, the element 
$\varphi +D^{n-1}\tilde \C^{n-1}$
is a coboundary in $\C^{n-1}$.

{\sc Case 2:} $\alpha =0$. Then \eqref{eq:Cocycle-n-1} implies
\begin{equation}\label{eq:5.3-prime}
\sum_{j=3}^{n}\sum_{k=1}^{j-2}(-1)^{k}i_j\alpha_{(i_1,\ldots,i_k+i_{k+1},\ldots,i_j-1,\ldots,i_n)}
+
\sum^{n-1}_{j=1}(-1)^{j+1}i_{j+1}\alpha_{(i_1,\ldots ,i_j+i_{j+1}-1,\ldots, i_n)}  = 0.
\end{equation}
For $i_n=1$, we obtain
\begin{align}
\nonumber
(-1)^n\alpha_{(i_1,i_2\ldots,i_{n-1})}&=\sum_{j=3}^{n-1}\sum_{k=1}^{j-2}(-1)^{k}i_j\alpha_{(i_1,\ldots,i_k+i_{k+1},\ldots,i_j-1,\ldots,i_{n-1},1)}\\
\nonumber
&+\sum_{k=1}^{n-2}(-1)^{k}i_j\alpha_{(i_1,\ldots,i_k+i_{k+1},\ldots,i_j-1,\ldots,i_{n-1},0)}\\
\label{eq:Cocycle-all-0,1}
 &+\sum^{n-2}_{j=1}(-1)^{j+1}i_{j+1}\alpha_{(i_1,\ldots ,i_j+i_{j+1}-1,\ldots ,i_{n-1},1)},
 \end{align}
and for $i_n=0$, $i_{n-1}=1$ it follows from \eqref{eq:5.3-prime} that
\begin{align}
\nonumber
  (-1)^{n-1}\alpha_{(i_1,i_2,\ldots,i_{n-2},0)}&=\sum_{j=3}^{n-1}\sum_{k=1}^{j-2}(-1)^{k}i_j\alpha_{(i_1,\ldots,i_k+i_{k+1},\ldots,i_j-1,\ldots,i_{n-2},1,0)}\\
  \label{eq:Cocycle-all-0}
     &+\sum^{n-2}_{j=1}(-1)^{j+1}i_{j+1}\alpha_{(i_1,\ldots ,i_j+i_{j+1}-1,\ldots ,i_{n-2},1,0)}.
\end{align}
Choose 
$\varphi_1\in \tilde \C^{n-2}$, 
$\psi_1\in \tilde \C^{n-1}$ in such a way that
\begin{equation}\label{eq:RedConstr}
\varphi_1[i_1|i_2|\ldots|i_{n-2}]=\beta_{(i_1,i_2,\ldots,i_{n-2})},
\quad 
\psi_1[i_1|i_2|\ldots|i_{n-1}]
=\sum\limits_{j=1}^{n-2}(-1)^{j}\beta_{(i_1,i_2,\ldots,i_j+i_{j+1},\ldots,i_{n-1})} ,
\end{equation}
where 
\[
\beta_{(i_1,i_2,\ldots,i_{n-3},0)}=
\begin{cases}
(-1)^{n-1}\alpha_{(i_1,i_2,\ldots,i_{n-3},1,0)},
   & i_1,i_2,\ldots,i_{n-3}\ge1,\ i_{n-2}=0, \\
0, & \text{otherwise}.
\end{cases}
\]
Then $\Delta^{n-2}\varphi_1-D^{n-1}\psi_1$ represents a coboundary in $C^{n-1}$, 
and
\[
\begin{aligned}
  (\Delta^{n-2}\varphi_1-D^{n-1}\psi_1)[i_1|i_2|\ldots |i_{n-1}]
  =
 &\sum_{j=3}^{n-1}\sum_{k=1}^{j-2}(-1)^{k}i_j\beta_{(i_1,\ldots,i_k+i_{k+1},\ldots,i_j-1,\ldots,i_{n-1})}\\
     &+\sum^{n-2}_{j=1}(-1)^{j+1}i_{j+1}\beta_{(i_1,\ldots ,i_j+i_{j+1}-1,\ldots ,i_{n-1})} .
\end{aligned}
\]
Hence,
\[
  (\Delta^{n-2}\varphi_1-D^{n-1}\psi_1)[i_1|i_2|\ldots|i_{n-3}|1|0]=\alpha_{(i_1,i_2,\ldots,i_{n-3},1,0)}=\varphi[i_1|i_2|\ldots|i_{n-3}|1|0].
\]
Therefore, we may assume that 
all $\alpha_{(i_1,i_2,\ldots,i_{n-3},1,0)}=0$. 
Relation \eqref{eq:Cocycle-all-0} implies 
$\alpha_{(i_1,i_2,\ldots,i_{n-2},0)}=0$  
for all $i_1,i_2,\ldots, i_{n-2}\ge 1$.

Let us repeat the construction \eqref{eq:RedConstr} 
with new values of $\beta$'s to get  
$\varphi'_1$ and $\psi'_1$ for 
\[
\beta{(i_1,i_2,\ldots,i_{n-2})}=
\begin{cases}
(-1)^n\alpha_{(i_1,i_2,\ldots,i_{n-2},1)},& i_1,i_2,\ldots,i_{n-2}\ge1, \\
0,& \text{otherwise}.
\end{cases}
\]
Then 
for all $i_1,\ldots,i_{n-2}\ge 1$ we get
\[
 (\Delta^{n-2}\varphi'_1-D^{n-1}\psi'_1)
   [i_1|i_2|\ldots|i_{n-3}|1|0]
   =0
\] 
and 
\[
 (\Delta^{n-2}\varphi'_1-D^{n-1}\psi'_1)
   [i_1|i_2|\ldots|i_{n-2}|1]
   =\alpha_{i_1,i_2,\ldots,i_{n-2},1}.
\]
Therefore, the cocycle $\varphi + D^{n-1}\tilde \C^{n-1}$
is a coboundary.
\end{proof}

\begin{corollary}\label{cor:AllFinite}
Let $M$ be a finite module over Weyl associative conformal algebra $U(2)$. Then $\Homol^k(U(2),M)=0$ for all $k\ge 2$.
\end{corollary}

\begin{proof}
Let $M$ be a finite conformal module over $U(2)$.
Then, in particular, $M$ is a finite module over the Virasoro Lie conformal algebra $\Vir $. 
Hence there exists a chain of $\Vir$-submodules
 (see, e.g., \cite[Lemma 3.3]{kolesnikov})
 \[
 0= M_{-1} \subset M_0 \subset \ldots \subset M_n=M,
 \]
where $M_i/M_{i-1}$,  $i=0,\ldots,n$, is 
either isomorphic to a $\Vir $-module $M_{(\alpha,\Delta)}$,
or trivial torsion-free module $\Bbbk [\partial ]u$ with 
$(v\oo\lambda u)=0$, 
or coincides with its torsion (hence, trivial).
Note that a $\Vir$-submodule of an $U(2)$-module $M$ is itself a 
$U(2)$-module, therefore, all $M_i$ are $U(2)$-modules and 
so are $M_i/M_{i-1}$. Hence, all irreducible quotients 
are of type $M_{(\alpha, 1)}$.

The case of torsion module was considered \cite{AlKol-JMP}. 
There is no difference between the scalar $U(2)$-module 
$\Bbbk $ 
and the trivial torsion-free module $M=\Bbbk [\partial ]u$. 
One may also apply the technique 
of Theorem~\ref{the5.5} to the case of trivial module $M$ as above 
to prove that $\Homol^k(U(2),M)=0$ for all $k\ge 2$.
Finally, Theorem~\ref{the5.5} implies
that
 \[
\Homol^k(U(2), M_i/M_{i-1} )=0, \quad i=0,\dots, n, 
\]
for all $k\ge 2$.
The short exact sequence of modules 
\[
0\to M_{i-1} \to M_i \to 
M_i/M_{i-1} \to 0
\]
leads to the long exact sequence of cohomology groups
\[
\begin{aligned}
\dots & \to \Homol^k(U(2),M_{i-1})
        \to \Homol^k(U(2),M_i)
        \to \Homol^k (U(2), {M_i/M_{i-1}}) \\
& \to \Homol^{k+1}(U(2),M_{i-1}) 
  \to \Homol^{k+1}(U(2),M_i)
        \to \Homol^{k+1} (U(2), {M_i/M_{i-1}}) \\
&        \to \dots .
\end{aligned}
\]
for every $i=1,\ldots,n$. 
Since  
\[
\Homol^k(U(2),M_0)=0,\quad \Homol^k(U(2),{M_1/M_0})=0
\]
for all  $k\ge2$, we obtain $\Homol^k(U(2), M_1)=0$.
Proceed by induction on $i=1,\dots, n$ to obtain
$\Homol^k(U(2),M_n)=0$, for all $k\ge2$.
\end{proof}

The series of modules $M_{(\alpha, \Delta)}$ also includes 
non-irreducible modules with $\Delta = 0$. 
For example, the module $M_{(0,0)}=\Bbbk [\partial ]u$
is defined by $(v\oo\lambda u) = \partial u$, 
this representation corresponds to the embedding 
of $\Vir $ into $\Cend_{x,1}$ \cite{BKL}.

\begin{corollary}
For all $n\ge 2$ we have     
$\Homol^n (U(2),M_{(\alpha, 0)}) = 0$.
\end{corollary}

Indeed, $M_{(\alpha,0)}$ is a finite module. The chain of submodules 
mentioned in the proof of Corollary~\ref{cor:AllFinite}
is given by 
\[
0=M_{-1} \subset (\partial+\alpha)M=M_0 \subset M=M_1, 
\]
where 
$M_0/M_{-1}\simeq M_{(\alpha, 1)}$, $M_1/M_0\simeq \Bbbk $.

\end{document}